\documentclass[10 pt]{amsart}
\usepackage{amsmath}
\usepackage{amssymb}
\usepackage{amsfonts}

\theoremstyle{plain}
\newtheorem{theorem}  {Theorem} []
\newtheorem{lemma}    [theorem] {Lemma} 
\newtheorem{corollary}[theorem] {Corollary}

\theoremstyle{definition}
\newtheorem{definition}[theorem]{Definition}

\theoremstyle{remark}
\newtheorem{remark}    [theorem]{Remark}
\newtheorem{example}   [theorem]{\bf Example}

\def \tdots {...}
\def \ttdots {\!...}

\begin{document}

\title
{Symmetrization of  Brace Algebras}
\author{Marilyn Daily}
\address{Max-Planck-Institute for Gravitational Physics, 
D-14476 Golm, Germany}
\email{Marilyn.Daily@aei.mpg.de}
\author{Tom Lada}
\address{Department of Mathematics, North Carolina State University,
Raleigh NC 27695}
\email{lada@math.ncsu.edu}
\thanks{This research was supported in part by NSF 
grant INT-0203119 and by grant M{\v S}MT ME 603}

\begin{abstract}
We show that the symmetrization of a brace algebra structure yields the structure of a symmetric brace algebra.  We also show that the symmetrization of the natural brace structure on $\bigoplus_{k\geq 1} Hom(V^{\otimes k},V)$ coincides with the natural symmetric brace structure on 
$\bigoplus_{k\geq 1} Hom(V^{\otimes k},V)^{as}$, the space of antisymmetric maps $V^{\otimes K}\rightarrow V$.
\end{abstract} 

\maketitle

\section{Introduction}

Brace algebras were first studied in the context of multilinear operations on the Hochschild complex of an associative algebra \cite{K88, G93, GV95}.  Symmetric brace algebras, in which the brace operations possess the property of graded symmetry, were subsequently introduced in \cite{LM04}.  Just as one may construct $L_{\infty}$ algebra structures by anti (skew) symmetrizing $A_{\infty}$ algebra structures \cite{LM95}, we show in this note that the symmetrization of a brace algebra structure yields a symmetric brace algebra structure.  
We prove in Section \ref{secTheorem1} that

$$
f\big< g_1,\dots,g_n\big> := \sum_{\sigma\in S_n} \epsilon(\sigma)
f\{g_{\sigma(1)},\dots,g_{\sigma(n)}\}
$$
where $\langle\:,\:\rangle$ and $\{\:,\:\}$ denote symmetric and non symmetric braces respectively.

The motivating example of a brace algebra is 
$\bigoplus_{k\geq 1} Hom(V^{\otimes k},V)$, and 
the fundamental example of a symmetric brace algebra is
the subspace of anti symmetric maps, 
$\bigoplus_{k\geq 1} Hom(V^{\otimes k},V)^{as}$.

In Section \ref{secTheorem2},
we show that these algebras 
are related by
$$
\sum_{\sigma\in S_n} \epsilon(\sigma) as(f\{g_{\sigma(1)},\dots,g_{\sigma(n)}\})
=
as(f)\big<as(g_1),\dots,as(g_n)\big>,
$$
where
$as(f)(v_1,\dots,v_k):=\sum_{\sigma\in S_k}
(-1)^{\sigma}\epsilon(\sigma)f(v_{\sigma(1)},\dots,v_{\sigma(k)})$ 
and $\epsilon(\sigma)$ is just the Koszul sign of the permutation.

In Sections \ref{secBrace} and \ref{secSymBrace}, 
we review the definitions and fundamental examples of brace algebras and symmetric brace algebras respectively.  
Section \ref{secLemmas} contains a collection of 
technical lemmas that are needed to prove the main theorems in the final two sections.

\medskip

\section{Brace Algebras}\label{secBrace}

\begin{definition}\label{defbrace1}
A {\em brace structure} on a graded vector space
consists of a collection of degree 0 multilinear braces
$x,x_1,\dots,x_n \mapsto x\{x_1,\dots x_n\}$ which satisfy
the identity, $x\{\,\} = x$, and in which 
$
x\{x_1,\!\tdots,x_n\}\{y_1,\!\tdots,y_r\}
$
is equal to
$$
\sum \!\epsilon \cdot 
x \{
y_1,\!\tdots,y_{i_1},
x_1\{y_{i_1+1},\!\tdots,y_{j_1}\},
y_{j_1+1},\!\tdots,y_{i_n},
x_n\{y_{i_n+1},\!\tdots,y_{j_n}\},
y_{j_n+1},\!\tdots,y_{r}
\}.
$$
In the above formula, the sum is over all sequences 
$0\!\leq\! i_1\!\leq\! j_1\!\leq\!\dots\!\leq\! i_n\!\leq\! j_n\!\leq\! r$, 
and $\epsilon$ is the Koszul sign of the permutation which maps
$(x_1,\ttdots,x_n,y_1,\ttdots,y_r)$ to
$$
(y_1,\ttdots,y_{i_1},x_1,y_{i_1+1},\ttdots,y_{j_1},
y_{j_1+1},\ttdots,y_{i_n},x_n,y_{i_n+1},\ttdots,y_{j_n},
y_{j_n+1},\ttdots,y_{r}).
$$
\end{definition}

The motivating example for a brace algebra structure is the space $Hom(V^{\otimes N},V)$ with the natural brace operation of degree $-n$ given by the composition
$$
f\{g_1,\dots,g_n\} 
= 
\!\!\!\!\!\!\!\!\!\!\!\!\!\!\!\!\!\!\!\!\!\!\!\!\!\!\!\!
\sum_{\qquad\qquad\ \ k_0+\dots+k_n=N-n}
\!\!\!\!\!\!\!\!\!\!\!\!\!\!\!\!\!\!\!\!\!\!\!\!\!\!\!\!
f\left(1^{\otimes k_0}\otimes g_1\otimes 1^{\otimes k_1}
\otimes\dots\otimes g_n\otimes 1^{\otimes k_n}\right),
$$
where $f\in Hom(V^{\otimes N},V)$.  This operation arises from the endomorphism
operad of $V$ considered in \cite{GV95}.  This operation was also utilized in
the context of the Hochschild complex of the associative algebra $V$ in
\cite{K88} and \cite{G93}.  After a regrading, this example may be regarded as a special case of the following

\begin{example}\label{exbrace}
Let $V$ be a graded vector space and consider the graded vector space $B_*(V)$ where 
$$B_s(V):=\bigoplus_{p-k+1=s}Hom(V^{\otimes k},V)_p$$
and where $Hom(V^{\otimes k},V)_p$ denotes the space of $k$-multilinear maps of degree $p$.
Given $f\in Hom(V^{\otimes N},V)_p$ and $g_i\in Hom(V^{\otimes a_i},V)_{q_i}$, define 
$f\{g_1,\dots,g_n\}\in Hom(V^{\otimes r},V)_{p+q_1+\dots+q_n}$ where $r=a_1+\dots+a_n+N-n$ by

$$
f\{g_1,\ttdots,g_n\} = 
\!\!\!\!\!\!\!\!\!\!\!\!
\sum_{k_0 + \dots +k_n=N-n} 
\!\!\!\!\!\!\!\!\!\!\!
(-1)^{\beta}
f(1^{\otimes k_0}\otimes g_1\otimes 1^{\otimes k_1}\otimes\dots\otimes
1^{\otimes k_{n-1}}\otimes g_n\otimes 1^{\otimes k_n}),
$$
where
$$
\beta=\sum_{j<i} \left[a_i-1\right]\left[k_j +  a_j\right]
+ \sum_{i} (N-i)\, q_i + \sum_{j<i} q_i\, a_j.
$$
\end{example}

\begin{remark}\label{Ainf}
In Example \ref{exbrace}, suppose that there exists a collection of maps
$$\mu_k\in Hom(V^{\otimes k},V)_{k-2}\in B_{-1}(V).$$
If we let $\mu=\mu_1+\mu_2+\dots$, then an $A_{\infty}$ algebra structure on $V$ may be described by the brace relation $\mu\{\mu\}=0$ \cite{LM04}.
\end{remark}

\medskip

\section{Symmetric Brace Algebras}\label{secSymBrace}

\begin{definition}
An {\em n-unshuffle} of $N$ elements is a partition 
$\sum_{i=1}^n a_i = N$ and a permutation $\gamma\in S_N$ such that
$$
\gamma(1) < \dots < \gamma(a_1), 
\gamma(1 + a_1) < \dots < \gamma(a_2+a_1),
\dots ,
\gamma\left(\!1+\sum_{i=1}^{n-1}\!a_i\!\right) < \dots < 
\gamma(N).
$$
\end{definition}

\begin{definition}\label{defsymbrace1}
A {\em symmetric brace algebra} is a graded vector space together with a
collection of degree zero multilinear braces 
$f\big<g_1,\dots,g_n\big>$ which are graded symmetric in $g_1,\dots,g_n$.
In a symmetric brace algebra, it is also required that $f\big<\big>=f$, 
and that $f\big<g_1,\dots,g_n\big>\big<x_1,\dots,x_r\big>$
be equal to
$$
\!\!\!\!\!\!\!\!\!\!\!\!\!
\sum_{\substack{\qquad\ \gamma \text{ is } (n+1) \\ \qquad\text{ unshuffle}}}\!\!\!\!\!\!\!\!\!\!\!\!\!\!
\epsilon \cdot f\big<g_1\big
<x_{\!\gamma(1)},\!\tdots,x_{\!\gamma(a_1)}\!\big>,\!\tdots,
g_n\big<x_{\!\gamma\left(\!1+\sum_{i=1}^{n-1}\!a_i\!\right)},
\!\tdots,x_{\!\gamma\left(\!\sum_{i=1}^{n}a_i\!\right)}\big>,
\!
x_{\!\gamma\left(\!1+\sum_{i=1}^{n}a_i\!\right)},\!\tdots,x_{\!\gamma(r)}\!\big>\!,
$$
where $\epsilon$ is the Koszul sign of the permutation which maps 
$(g_1,\tdots,g_n,x_1,\tdots,x_r)$
to
$$
\left(g_1,x_{\gamma(1)},\dots,
x_{\!\gamma(a_1)},g_2,\!\tdots,
x_{\!\gamma\left(\!1+\sum_{i=1}^{n-1}\!a_i\!\right)},
\!\tdots,x_{\!\gamma\left(\!\sum_{i=1}^{n}a_i\!\right)},g_n,
\!
x_{\!\gamma\left(\!1+\sum_{i=1}^{n}a_i\!\right)},\!\tdots,x_{\!\gamma(r)}\right).
$$
\end{definition}

Just as with brace algebras, the fundamental example of a symmetric brace algebra is provided by the space of antisymmetric maps of degree $p$, $Hom(V^{\otimes k},V)^{as}_p$.  To be precise, we have

\begin{example}\label{exsymbrace}
Let $V$ be a graded vector space and $B_*(V)$ be the graded vector space given by 
$$B_s(V)=\bigoplus_{p-k+1=s}Hom(V^{\otimes k},V)^{as}_p,$$
Given $f \in Hom(V^{\otimes k},V)^{as}_p$ and $g_i \in Hom(V^{\otimes
a_i},V)^{as}_{q_i}$, $1 \leq i \leq n$,
define the symmetric brace
$$
f\big<g_1,\dots,g_n\big>(x_1,\dots,x_r) 
= (-1)^{\delta}
\!\!\!\!\!\!\!\!\!\!\!\!
\sum_{\substack{\gamma \text{ is an }\\ 
(a_1\vert a_2\vert\dots\vert a_{n+1}) \\ 
\text{ unshuffle}}}\!\!\!\!\!\!\!\!\!\!\!\!
\chi(\gamma)
f(g_1\otimes\dots\otimes g_n\otimes 1^{\otimes N-n})
(x_{\gamma(1)},\dots,x_{\gamma(r)}),
$$
where 
$$
\delta=\sum_i ^n (N-i) q_i +\sum_{j<i}q_i a_j + \sum_{j<i}a_i a_j + \sum_i (n-i)a_i,
$$

and $\chi$ is the antisymmetric Koszul 
sign of the permutation $\gamma$.
\end{example}

\begin{remark}\label{Linf}
Suppose that in Example \ref{exsymbrace} we have maps
$$ l_k\in Hom(V^{\otimes k},V)^{as}_{k-2}\in B_{-1}(V).$$
If we let $l=l_1+l_2+\dots$, then an $L_{\infty}$ algebra structure on $V$ is given by the symmetric brace relation $l\langle l\rangle=0$.
\end{remark}

\medskip

\section{Some Lemmas}\label{secLemmas}

Although the expressions in this paper involve many sums, 
permutations, and antisymmetrizations, we will be able to 
simplify things considerably with the help of the following lemmas.
Lemma \eqref{techlem1} provides a decomposition of $as(f)$
which will be useful later.

\begin{lemma}\label{techlem1}
$as(f)=f\circ\Phi_{nm}\circ\Psi_n\circ\Theta_m\ \
\forall\ f\in Hom(V^{\otimes n+m},V)$, where
\begin{align*}
\Theta_m (y_1,\ttdots,y_n, z_1,\ttdots z_m)
&=
\sum_{\pi\in S_m}\chi(\pi) (y_1,\ttdots,y_n, z_{\pi(1)},\ttdots z_{\pi(m)}),
\\
\Psi_n (y_1,\ttdots,y_n, z_1,\ttdots z_m)
&=
\sum_{\sigma\in S_n}\chi(\sigma)
(y_{\sigma(1)},\ttdots,y_{\sigma(n)}, z_1,\ttdots z_m),
\\
\Phi_{\!nm\!} (y_1,\ttdots,y_n, z_1,\ttdots, z_m) 
&= 
\!\!\!\!\!\!\!\!\!\!\!\!\!\!\!\!\!\!\!\!\!\!
\sum_{\qquad\quad k_0+\dots +k_n = m} 
\!\!\!\!\!\!\!\!\!\!\!\!\!\!\!\!\!\!\!\!
(\!-\!1)^{\eta}\!
\left( 
z_{1},\ttdots, z_{k_0},
y_{1},
z_{1+k_0},
\ttdots,
y_{n},
z_{1+k_0 + \dots + k_{n-1}},\ttdots, z_{m}
\right )\!\!,
\end{align*}
and
$
\eta = \sum_{i=1}^n \{ y_{i}\left[
z_{1} + \dots + z_{(k_0 + k_1 + \dots + k_{i-1})}\right] + 
(n-i) k_i \}.
$
\end{lemma}

\begin{proof}

Since $\Psi_n$ does all permutations of the first $n$ inputs, 
$\Theta_m$ provides all permutations of the last $m$ inputs, and
$\Phi_{nm}$ distributes the last $n$ variables between
the first $m$ in every possible way, the composition is clearly a
sum of all permutations of the original $n+m$ variables.
A moment's reflection also reveals that the sign of each summand
in the composition
is the Koszul sign together with the sign of the permutation.
\end{proof}

Lemma \eqref{techlem2} states that if we sum over all (signed)
$(a_1\vert\dots\vert a_n)$ unshuffles, and then sum over all (signed)
permutations of the $a_i$ variables in each piece, then this is equivalent
to just summing over all signed permutations of the original 
$a_1 + \dots + a_n$ variables.

\begin{lemma}\label{techlem2}
If $N=a_1+\dots + a_{n}$,
then
$\sum_{\pi\in S_{N}}\chi(\pi)\!\left(\!x_{\pi(1)},\dots, x_{\pi(N)}\right)$
is equal to
$$
\!\!\!\!\!\!\!\!\!\!\!\!\!\!
\sum_{\substack{\ \gamma \text{ is }\\ 
\qquad\ (a_1\vert\dots\vert a_{n}) \\ 
\qquad\! \text{ unshuffle}}}\!\!\!\!\!\!\!\!\!\!\!\!\!\!
\chi(\gamma)
\!\!\!\!\!\!\!\!\!\!
\sum_{\ \ \ \ \pi_1\in S_{a_1}}
\!\!\!\!\!\!\!\!
\chi(\pi_1)
\tdots
\!\!\!\!\!\!\!\!\!\!\!\!\!
\sum_{\qquad \pi_{a_{n}}\in S_{a_n}}
\!\!\!\!\!\!\!\!\!\!\!\!\!
\chi(\pi_{n})\!
\!\left(\!
x_{\gamma(\pi\!(1)\!)},\ttdots,\! x_{\gamma(\pi_1\!(a_1)\!)},
x_{\gamma(\pi_2(1)+a_1)},
\ttdots,
x_{\gamma\left(\pi_{\!n}\!(a_n)+\sum_{i=1}^{n\!-\!1} a_i\right)}
\!\right)\!.
$$
\end{lemma}

\begin{proof}
Clearly, the right hand side is the sum of distinct permutations of the
$x$ terms with the correct sign.  Furthermore, since there are 
$\tfrac{N!}{(a_1)!\dots (a_n)!}$ unshuffles $\gamma$ and $(a_i)!$
permutations $\pi_i$, there are $N!$ summands in the right hand side, 
which agrees with the number of summands on the left hand side.
\end{proof}

\begin{lemma}\label{TechLemPiHat}
Suppose $k_0+a_1+k_1+\dots+a_n+k_n=r$, 
$\sigma\in S_n$, and $\pi\in S_r$.  
Let $A=a_1+\dots a_n$, denote
$
X_i = x_{\pi\left(1+a_1 + \dots + a_{i-1}\right)},\dots, 
x_{\pi\left(a_1 + \dots + a_{i}\right)},
$
\phantom{\Big(}\!\!\!
and also denote 
$
X_{\pi} = 
x_{\!\pi(1+A)},\ttdots,x_{\!\pi(k_0+A)}, X_{\sigma(1)},
x_{\!\pi(1+k_0+A)},\ttdots,
X_{\sigma(n)},
x_{\!\pi(1+k_0+\dots+k_{n-1}+A)},\ttdots,x_{\!\pi(r)}.
$
Then
we can define $\hat{\pi}\in S_r$ by \phantom{\Big(}\!\!\!
$$
\hat{\pi}(i) = 
\left\{
\begin{aligned}
\pi\!\left(\!i+A-\sum_{j\leq m}a_{\sigma(j)}\!\right)
&\text{ if }
\sum_{j<m} k_j + \sum_{j\leq m} a_{\sigma(j)} 
< i \leq
\sum_{j\leq m} k_j + \sum_{j\leq m} a_{\sigma(j)}.
\\
\pi\!\left(\!i-\!\sum_{j<m} k_j + \!\sum_{j<\sigma(m)}\!a_{j}\!\right)
&\text{ if }
\sum_{j<m} k_j + \sum_{j<m} a_{\sigma(j)} 
< i \leq
\sum_{j<m} k_j + \sum_{j\leq m} a_{\sigma(j)}.
\end{aligned}
\right.
$$
Furthermore, given this notation, 
\begin{align*}
&\ \ \ 
X_{\pi} = x_{\hat{\pi}(1)},\dots,x_{\hat{\pi}(r)} 
\quad\text{and}\quad
\epsilon(\hat{\pi}) = \epsilon(\pi)(-1)^{\alpha_1}
\quad\text{and}\quad
\chi(\hat{\pi}) = \chi(\pi)(-1)^{\alpha_2},
\\
&\text{where }
\alpha_1 = 
\!\!\!\!\!\!\!\!\!\!\!\!\!\!\!\!\!\!\!\!\!\!\!\!\!\!
\sum_{\qquad\qquad\ 
i<j\ \&\ \sigma(i) > \sigma(j)}
\!\!\!\!\!\!\!\!\!\!\!\!\!\!\!\!\!\!\!\!\!\!\!\!\!
\vert X_{\sigma(i)}\vert\ \vert X_{\sigma(j)}\vert 
\ +\,
\sum_{i=1}^n 
\vert X_{\sigma(i)}\vert 
\left[x_{\pi(1+A)}+ \dots + x_{\pi(k_0 + \dots + k_{i-1}+A)}\right]
\\
&\text{and }\ 
\alpha_2 = \alpha_1 + 
\!\!\!\!\!\!\!\!\!\!\!\!\!\!\!\!\!\!\!\!\!\!\!\!\!\!
\sum_{\qquad\qquad\ 
i<j\ \&\ \sigma(i) > \sigma(j)}
\!\!\!\!\!\!\!\!\!\!\!\!\!\!\!\!\!\!\!\!\!\!\!\!\!
a_{\sigma(i)} a_{\sigma(j)}
\ +\,
\sum_{j<i} a_{\sigma(i)} k_j.
\end{align*}
\end{lemma}

\begin{proof}
Careful examination of the definition of $\hat{\pi}$ reveals that
the first formula
moves \lq\lq free" strings
of the form $x_{\pi(1+k_0+\dots+k_{i-1})},\dots,x_{\pi(k_0+\dots+k_{i})}$
into place (for $0\leq m\leq n$), 
and the second formula 
relocates the strings $X_{\sigma(i)}$
(for $1\leq m\leq n$).
Thus $X_{\pi} = x_{\hat{\pi}(1)},\dots,x_{\hat{\pi}(r)}$.

Furthermore, when $x_{\pi(1)},\dots,x_{\pi(r)}$ are permuted to yield
$x_{\hat{\pi}(1)},\dots,x_{\hat{\pi}(r)}$, the Koszul sign is 
$(-1)^{\alpha_1}$,  where the first sum in $\alpha_1$ comes from 
$\sigma$ 
permuting the $X_i$ strings, and the second sum comes from moving the
\lq\lq free" strings into place.  Finally, the additional sums in 
$\alpha_2$ count the transpositions, yielding the correct 
antisymmetric Koszul sign.
\end{proof}

\begin{lemma}\label{techlem3}
Suppose that $\sigma\in S_n$ permutes $\{v_1\dots v_n\}$ and $\{w_1\dots w_n\}$.
Then
\begin{align*}
&(1)\
\sum_{i>j} v_i w_j + 
\!\!\!\!\!\!\!\!\!\!\!\!\!\!\!\!\!\!\!\!\!\!
\sum_{\qquad\qquad\ 
i<j\ \&\ \sigma(i) > \sigma(j)}
\!\!\!\!\!\!\!\!\!\!\!\!\!\!\!\!\!\!\!\!\!\!\!\!\!
\left\{ w_{\sigma(i)} v_{\sigma(j)} + v_{\sigma(i)} w_{\sigma(j)}\right\}
+ \sum_{i>j} v_{\sigma(i)} w_{\sigma(j)} \equiv 0\ (mod\ 2).
\\
&(2)\
\!\!\!\!\!\!\!\!\!\!\!\!\!\!\!\!\!\!\!\!\!\!\!\!\!\!\!
\sum_{\qquad\qquad\ 
i<j\ \&\ \sigma(i) > \sigma(j)}
\!\!\!\!\!\!\!\!\!\!\!\!\!\!\!\!\!\!\!\!\!\!\!\!\!
\left\{ v_{\sigma(i)} + v_{\sigma(j)}\right\}
\equiv \sum_i(i-1) v_i + \sum_i (i-1) v_{\sigma(i)}\ (mod\ 2).
\end{align*}
\end{lemma}

\begin{proof}
To prove the first assertion, we note that 
$$
\!\!\!\!\!\!\!\!\!\!\!\!\!\!\!\!\!\!\!\!\!\!\!\!\!\!\!\!
\sum_{\substack{\qquad\qquad\ \ \ i < j\ \&\  \sigma(i) > \sigma(j)}}
\!\!\!\!\!\!\!\!\!\!\!\!\!\!\!\!\!\!\!\!\!\!\!\!\!\!\!\!\!\!
\left\{ v_{\sigma\!(\!i\!)} w_{\sigma\!(\!j\!)} 
+ w_{\sigma\!(\!i\!)} v_{\sigma\!(\!j\!)}\!\right\}
+ \sum_{i>j} v_{\sigma\!(\!i\!)} w_{\sigma\!(\!j\!)}
=
\!\!\!\!\!\!\!\!\!\!\!\!\!\!\!\!\!\!
\sum_{\substack{\qquad\quad  i \!<\! j \& \sigma\!(\!i\!) > \sigma\!(\!j\!)}}
\!\!\!\!\!\!\!\!\!\!\!\!\!\!\!\!\!\!
v_{\sigma\!(\!i\!)} w_{\sigma\!(\!j\!)}
+
\!\!\!\!\!\!\!\!\!\!\!\!\!\!\!\!\!\!
\sum_{\substack{\qquad\quad i \!>\! j \&  \sigma\!(\!i\!) < \sigma\!(\!j\!)}}
\!\!\!\!\!\!\!\!\!\!\!\!\!\!\!\!\!\!
v_{\sigma\!(\!i\!)} w_{\sigma\!(\!j\!)} 
+ 
\sum_{i>j} v_{\sigma\!(\!i\!)} w_{\sigma\!(\!j\!)},
$$
$$
\text{which is congruent } (mod\ 2) \text{ to} 
\!\!\!\!\!\!\!\!\!\!\!\!\!\!\!\!
\sum_{\substack{\qquad\quad i\!<\! j \& \sigma\!(\!i\!) > \sigma\!(\!j\!)}}
\!\!\!\!\!\!\!\!\!\!\!\!\!\!\!\!\!\!
v_{\sigma\!(\!i\!)} w_{\sigma\!(\!j\!)}
+
\!\!\!\!\!\!\!\!\!\!\!\!\!\!\!\!\!
\sum_{\substack{\qquad\quad i\! >\! j \&  \sigma\!(\!i\!) > \sigma\!(\!j\!)}}
\!\!\!\!\!\!\!\!\!\!\!\!\!\!\!\!\!\!\!
v_{\sigma\!(\!i\!)} w_{\sigma\!(\!j\!)}
=
\!\!\!\!\!\!\!
\sum_{\quad\sigma\!(\!i\!)>\sigma\!(\!j\!)}\!\!\!\!\!\!\!
v_{\sigma\!(\!i\!)} w_{\sigma\!(\!j\!)}
=
\sum_{i>j} v_{i} w_{j}.
$$
To prove the second statement, suppose that all $w_i$ are odd.
Then
$$
\!\!\!\!\!\!\!\!\!\!\!\!\!\!\!\!\!\!\!\!\!\!\!\!\!\!\!\!
\sum_{\substack{\qquad\qquad\ \ \ i < j\ \&\  \sigma(i) > \sigma(j)}}
\!\!\!\!\!\!\!\!\!\!\!\!\!\!\!\!\!\!\!\!\!\!\!\!\!\!\!\!
\left\{ v_{\sigma(i)} + v_{\sigma(j)} \right\}
\equiv
\!\!\!\!\!\!\!\!\!\!\!\!\!\!\!\!\!\!\!\!\!\!\!\!\!\!\!\!
\sum_{\substack{\qquad\qquad\ \ \ i < j\ \&\  \sigma(i) > \sigma(j)}}
\!\!\!\!\!\!\!\!\!\!\!\!\!\!\!\!\!\!\!\!\!\!\!\!\!\!\!\!
\left\{ v_{\sigma(i)} w_{\sigma(j)} + w_{\sigma(i)} v_{\sigma(j)}\right\}
\equiv
\sum_{i>j}
\left\{ v_{i} w_{j} + v_{\sigma(i)} w_{\sigma(j)}\right\}
$$
(by the first assertion).  Since all $w$-terms are odd, this is congruent to
$$
\sum_{j=1}^n\sum_{\ i=j+1}^n (v_i + v_{\sigma(i)})
=
\sum_j (j-1) v_j + \sum_j (j-1) v_{\sigma(j)}.
$$
\end{proof}

\medskip

\section{Symmetrization of Brace Algebras}\label{secTheorem1}

Given a (non-symmetric) brace structure $\{ ,\}$ on a graded vector space,
we can define a symmetric brace structure $\big<,\big>$ 
via
$$
f\big< g_1,\dots,g_n\big> := \sum_{\sigma\in S_n} \epsilon(\sigma)
f\{g_{\sigma(1)},\dots,g_{\sigma(n)}\}.
$$
Clearly, this satisfies the first symmetric brace axiom, since
$f\big<\,\big>=f\{\,\}=f$.  We show in Theorem \eqref{Theorem1} 
that it satisfies the second symmetric brace axiom given in 
Definition \eqref{defsymbrace1}, so this does in fact induce a
symmetric brace structure.  First, however, we need the following 
two lemmas, which are analogous to Lemmas \eqref{techlem1} and \eqref{techlem2}.

\begin{lemma}\label{techlem1a}\ 
$\sum_{\rho\in S_{n+m}} \epsilon(\rho)\, f\{x_{\rho(1)},\dots,x_{\rho(n)}\}
\,=\,
\tilde{f}_n\circ\theta_m(x_1,\dots,x_{n+m})$, 
where
\phantom{\Bigg(}
$
\Theta_m (y_1,\ttdots,y_n, z_1,\ttdots z_m) =
\sum_{\pi\in S_m}\epsilon(\pi) (y_1,\ttdots,y_n, z_{\pi(1)},\ttdots z_{\pi(m)})$
and
\begin{multline*}
\tilde{f}_{n} (y_1,\ttdots,y_n, z_1,\ttdots, z_m)
\\[2pt]=
\!\!\!
\sum_{\ \ \sigma\in S_n}
\!\!\!\!
\epsilon(\sigma)
\!\!\!\!\!\!\!\!\!\!\!\!\!\!\!\!\!\!\!\!
\sum_{\qquad\quad k_0+\dots +k_n = m}
\!\!\!\!\!\!\!\!\!\!\!\!\!\!\!\!\!\!\!\!
(-1)^{\eta}
f\!\left\{
z_{1},\ttdots, z_{k_0},
y_{\sigma(1)},
z_{1+k_0},
\ttdots,
y_{\sigma(n)},
z_{1+k_0 + \dots + k_{n-1}},\ttdots, z_{m}
\right\}\!\!,
\end{multline*}
with a Koszul sign given by 
$\eta = \sum_{i=1}^n  y_{\sigma(i)}\left[
z_{1} + \dots + z_{(k_0 + k_1 + \dots + k_{i-1})}\right]$.
\end{lemma}

\begin{lemma}\label{techlem2a}
If $N=a_1+\dots + a_{n}$,
then
$\sum_{\pi\in S_{N}}\epsilon(\pi)\!\left(\!x_{\pi(1)},\dots, x_{\pi(N)}\right)$
is equal to
$$
\!\!\!\!\!\!\!\!\!\!\!\!\!\!
\sum_{\substack{\ \gamma \text{ is }\\ 
\qquad\ (a_1\vert\dots\vert a_{n}) \\ 
\qquad\! \text{ unshuffle}}}\!\!\!\!\!\!\!\!\!\!\!\!\!\!
\epsilon(\gamma)
\!\!\!\!\!\!\!\!\!\!
\sum_{\ \ \ \ \pi_1\in S_{a_1}}
\!\!\!\!\!\!\!\!
\epsilon(\pi_1)
\tdots
\!\!\!\!\!\!\!\!\!\!\!\!\!
\sum_{\qquad \pi_{a_{n}}\in S_{a_n}}
\!\!\!\!\!\!\!\!\!\!\!\!\!
\epsilon(\pi_{n})\!
\!\left(\!
x_{\gamma(\pi\!(1)\!)},\ttdots,\! x_{\gamma(\pi_1\!(a_1)\!)},
x_{\gamma(\pi_2(1)+a_1)},
\ttdots,
x_{\gamma\left(\pi_{\!n}\!(a_n)+\sum_{i=1}^{n\!-\!1} a_i\right)}
\!\right)\!.
$$
\end{lemma}

\begin{remark}
Although a brace structure allows operators $g$ which accept an arbitrary
number of inputs, it will be convenient in the proof of the following theorem 
to let $g^a$ denote the restriction of $g$ which accepts only exactly $a$ inputs.
\end{remark}

\begin{theorem}\label{Theorem1}
Given a (non-symmetric) brace structure $\{\, ,\}$ on a graded vector space,
define $\big<\, ,\big>$ 
via
$$
f\big< g_1,\dots,g_n\big> := \sum_{\sigma\in S_n} \epsilon(\sigma)
f\{g_{\sigma(1)},\dots,g_{\sigma(n)}\}.
$$
Then
$f\big<g_1,\dots,g_n\big>\big<x_1,\dots,x_r\big>$
is equal to
$$
\!\!\!\!\!\!\!\!\!\!\!\!\!
\sum_{\substack{\qquad\ \gamma \text{ is } (n+1) \\ \qquad\text{ unshuffle}}}
\!\!\!\!\!\!\!\!\!\!\!\!\!\!
\epsilon \cdot f\big<g_1\big
<x_{\!\gamma(1)},\!\tdots,x_{\!\gamma(a_1)}\!\big>,\!\tdots,
g_n\big<x_{\!\gamma\left(\!1+\sum_{i=1}^{n-1}\!a_i\!\right)},
\!\tdots,x_{\!\gamma\left(\!\sum_{i=1}^{n}a_i\!\right)}\big>,
\!
x_{\!\gamma\left(\!1+\sum_{i=1}^{n}a_i\!\right)},\!\tdots,x_{\!\gamma(r)}\!\big>\!,
$$
where $\epsilon$ is the Koszul sign of the permutation which maps 
$(g_1,\tdots,g_n,x_1,\tdots,x_r)$
to
$$
\left(g_1,x_{\gamma(1)},\dots,
x_{\!\gamma(a_1)},g_2,\!\tdots,
x_{\!\gamma\left(\!1+\sum_{i=1}^{n-1}\!a_i\!\right)},
\!\tdots,x_{\!\gamma\left(\!\sum_{i=1}^{n}a_i\!\right)},g_n,
\!
x_{\!\gamma\left(\!1+\sum_{i=1}^{n}a_i\!\right)},\!\tdots,x_{\!\gamma(r)}\right).
$$
\end{theorem}

\begin{proof}
First, we will look at the right hand side.

If we temporarily denote 

\begin{align*}
h_k &= g_k\big<x_{\gamma\left(1+a_1 + \dots + a_{k-1}\right)},\ttdots, 
x_{\gamma\left(a_1 + \dots + a_{k}\right)}\big>
\\[2pt]
&= 
\sum_{\ \pi_k\in S_{a_k}}\!\!\!\epsilon(\pi_k)\
g_k\left\{ x_{\gamma\left(\pi_k(1)+a_1 + \dots + a_{k-1}\right)},\ttdots, 
x_{\gamma\left(\pi_k(a_k) + a_1 + \dots + a_{k-1}\right)}\!\right\},
\end{align*}

and denote $A=\sum_{i=1}^{n} a_i$,
then the right hand side is equal to
$$
\!\!\!\!\!\!\!\!\!\!\!\!\!\!\!\!\!\!\!\!\!\!\!\!\!\!\!\!\!\!\!\!\!\!\!\!\!\!
\sum_{\substack{\qquad\qquad\quad
 a_1 + \dots + a_{n+1} = r\ \&\ 
 \\
 \qquad\qquad\qquad\quad\
 \gamma \text{ is }  
(a_1\vert\dots\vert a_{n+1})
\text{ unshuffle}}}
\!\!\!\!\!\!\!\!\!\!\!\!\!\!\!\!\!\!\!\!\!\!\!\!\!\!\!\!\!\!\!\!\!\!\!\!\!
\!\!\!\!\!\!\!
(\!-\!1)^{\nu}\ \epsilon(\gamma)\ f\big<
h_1,\dots,h_n,
x_{\!\gamma\left(\!1+A\!\right)},\tdots,x_{\!\gamma(a_{n+1}+A)}\big>,
$$
where
$\nu = \sum_{i=2}^n g_i [x_{\gamma(1)}+\dots + x_{\gamma(a_1 +\dots +a_{i-1})}]$
is a Koszul sign.
\phantom{\Big(}
After applying 
Lemma \eqref{techlem1a}, this is equal to 
$$
\!\!\!\!\!\!\!\!\!\!\!\!\!\!\!\!\!\!\!\!\!\!\!\!\!
\sum_{\substack{\qquad\qquad\ 
 a_1 + \dots + a_{n+1} = r,
 \\
 \qquad\quad\ 
 \gamma \text{ is unshuffle}}}
\!\!\!\!\!\!\!\!\!\!\!\!\!\!\!\!\!\!\!\!\!\!\!
(\!-\!1)^{\nu}\ \epsilon(\gamma)\
\tilde{f}_n\left(
\!\!\!\!\!\!\!\!\!\!\!\!\!\!\!\!\!\!\!\!\!\!\!\!
\sum_{\qquad\qquad\pi_{n+1}\in S_{a_{n+1}}}
\!\!\!\!\!\!\!\!\!\!\!\!\!\!\!\!\!\!\!\!\!\!\!
\epsilon(\pi_{n+1}) 
(h_{1},\ttdots, h_{n},
x_{\gamma(\pi_{n+1}(1)+A)},\ttdots,x_{\gamma(\pi_{n+1}(a_n+1)+A)}
\!\right)\!\!,
$$
where $\tilde{f_n}$ is as defined in Lemma \eqref{techlem1a}.
Now, we will pull all of the $x$ terms back out, in order to apply 
Lemma \eqref{techlem2a}.  Note that the Koszul signs from this 
transformation merely cancel out $(-1)^{\nu}$.
We then have the following long formula:
\vspace{4pt}

\begin{multline*}
\!\!\!\!\!\!\!\!\!\!
\sum_{\ \ 
 (a_i),\gamma}\!\!\!\!\!
\epsilon(\gamma)\
\!\!\!\!\!\!\!\!\!\!
\sum_{\quad\ \pi_{1}\in S_{a_1}}
\!\!\!\!\!\!\!\!\!
\epsilon(\pi_{1}) 
\tdots\!\!\!\!\!\!\!\!\!\!\!\!\!\!\!\!\!\!\!\!\!\!\!\!\!
\sum_{\qquad\qquad \pi_{(n+1)}\in S_{a_{n+1}}}
\!\!\!\!\!\!\!\!\!\!\!\!\!\!\!\!\!\!\!\!\!\!\!\!
\epsilon(\pi_{n+1}) 
\tilde{f}_n\!\left(
g_1^{a_1},\ttdots,g_n^{a_n}\!,\!1^{a_{n+1}}
\right)\!
\big(\!
x_{\!\gamma\left(\pi_{\!1}(1)\!\right)},\ttdots, 
x_{\!\gamma\left(\pi_{\!1}(a_1)\!\right)},
\\[-8pt]
\phantom{************************,}
x_{\!\gamma\left(\pi_2(1)+a_1\!\right)\!},
\ttdots, 
x_{\!\gamma(\pi_{n\!}(A\!)},
\\
\phantom{*****************************,}
x_{\!\gamma(\!\pi_{\!n\!+\!1}(1)+A\!)},\ttdots
x_{\!\gamma(\pi_{n\!+\!1}(a_n+1)+A\!)\!}
\big)\!.
\end{multline*}
Now, though, we can apply Lemma \eqref{techlem2a},
which yields the much shorter formula,
$$
\!\!\!\!\!\!\!\!\!\!
\sum_{\ 
 (a_i)}
\!\!\!
\sum_{\quad \pi\in S_{r}}
\!\!\!\!
\epsilon(\pi)\,
\tilde{f}_n\!\left(
g_1^{a_1},\ttdots,g_n^{a_n}\!,\!1^{a_{n+1}}
\right)\!
\left(\!
x_{\!\pi(1)},\ttdots, 
x_{\!\pi(r)}
\right)\!.
$$
Before continuing, we need to pull all of the $x$ terms back inside.
In order to make our expressions a bit shorter,
let $X_i$ denote the input to $g_i$.  In other words, define
$$
X_i = x_{\pi\left(1+a_1 + \dots + a_{i-1}\right)},\dots, 
x_{\pi\left(a_1 + \dots + a_{i}\right)}\quad for\ i\in\{1\dots n\}.
$$
It will also be convenient to let $\vert X_i\vert$ denote the sum
of the degrees of the variables in $X_i$.
When we pull the $x$-terms inside and use the more concise notation 
just defined, the formula for the right hand side becomes
$$
\!\!\!\!\!\!\!\!\!\!
\sum_{\ 
 (a_i)}
\!\!
\sum_{\ \ \pi\in S_{r}}
\!\!\!
\epsilon(\pi)\,
(-1)^{\tilde{\nu}}
\tilde{f}_n\!\left(
g_1(X_1),\ttdots,g_n(X_n),
x_{\pi(1+A)},\dots,x_{\pi(r)}
\right)\!,
$$
where $\tilde{\nu} = \sum_{j<i} g_i \vert X_j\vert$.
After expanding $\tilde{f}_n$,
the right hand side is equal to
\vspace{2pt}
\begin{multline*}
\!\!\!\!\!\!\!
\sum_{
 (a_i)}
\!\!\!
\sum_{\ \ \pi\in S_{r}}
\!\!\!\!\!
\epsilon(\pi)
(\!-\!1)^{\tilde{\nu}\!}
\!\!\!\!
\sum_{\ \ \sigma\in S_n}
\!\!\!\!
\epsilon(\sigma\!)
\!\!\!\!\!\!\!\!\!\!\!\!\!\!\!\!\!\!\!\!\!\!\!\!\!\!
\sum_{\qquad\qquad k_0+\dots +k_n = a_{n+1}}
\!\!\!\!\!\!\!\!\!\!\!\!\!\!\!\!\!\!\!\!\!\!\!\!\!
(\!-\!1)^{\eta}\!
f\big\{
x_{\pi(1+A)},\ttdots,x_{\pi(k_0+A)},
g_{\sigma(1)}(X_{\sigma(1)}),
x_{\pi(1+k_0+A)},
\ttdots\\[-6pt]
\phantom{************************}
\ttdots,
g_{\sigma(n)}(X_{\sigma(n)}),
x_{\pi(1+ k_0 + \dots + k_{n}+A)},
\ttdots,
x_{\pi(r)}
\big\}\!.\!\!\!\!\!
\\[3pt]
\text{Here, }
\eta = 
\sum_{i=1}^n
\left( g_{\sigma(i)} + \vert X_{\sigma(i)}\vert\right)
\left[ x_{\pi(1+A)}+\dots+x_{\pi(k_0+\dots+k_{i-1}+A)}\right]
\phantom{*********,}
\\
\text{and }
\epsilon(\sigma) = (-1)^{\lambda},\ \text{where }
\lambda
=
\!\!\!\!\!\!\!\!\!\!\!\!\!\!\!\!\!\!\!\!\!\!\!\!\!\!\!\!\!
\sum_{\qquad\qquad\ \ \ i<j\  \&\ \sigma(i) > \sigma(j)} 
\!\!\!\!\!\!\!\!\!\!\!\!\!\!\!\!\!\!\!\!\!\!\!\!\!\!
\left( g_{\sigma(i)} + \vert X_{\sigma(i)}\vert\right)
\left( g_{\sigma(j)} + \vert X_{\sigma(j)}\vert\right).
\phantom{*******}
\end{multline*}

Now, we will look at the left hand side.
$f\big<g_1,\dots,g_n\big>\big<x_1,\dots,x_r\big>$
is equal to
$\sum_{\sigma\in S_n} \epsilon(\sigma) f\{g_{\sigma(1)},\dots,g_{\sigma(n)}\}
\big<x_1,\dots,x_r\big>$,
which is equal to
$$
\sum_{\sigma\in S_n} \epsilon(\sigma) \sum_{\pi\in S_r} \epsilon(\pi)\
f\{g_{\sigma(1)},\dots,g_{\sigma(n)}\}
\{x_{\pi(1)},\dots,x_{\pi(r)}\}.
$$
If we apply Definition \eqref{defbrace1} and let $g_i^{a_i}$ denote 
the restriction of $g_i$ which accepts exactly $a_i$ inputs,
then the left hand side is equal to
$$
\!\!\!\!\!\!\!\!
\sum_{\ \ \sigma\in S_n} \!\!\!\! \epsilon(\sigma) \!\!\!
\sum_{\ \pi\in S_r} \!\! \epsilon(\pi)
\!\!\!\!\!\!\!\!\!\!\!\!\!\!\!\!\!\!\!\!\!\!\!
\!\!\!\!\!\!\!\!\!\!\!\!\!\!\!\!\!\!\!\!\!\!
\sum_{\qquad \qquad\qquad\quad k_0+\dots + k_n + a_1 + \dots + a_n = r} 
\!\!\!\!\!\!\!\!\!\!\!\!\!\!\!\!\!\!\!\!\!\!
\!\!\!\!\!\!\!\!\!\!\!\!\!\!\!\!\!\!\!\!\!
f\{1^{k_0},g_{\sigma(1)}^{a_\sigma(1)},1^{k_1},\dots,g_{\sigma(n)}^{a_\sigma(n)},1^{k_n}\}
(x_{\pi(1)},\dots,x_{\pi(r)}).
$$
After applying Lemma \eqref{TechLemPiHat}, this is equal to

\begin{multline*}
\!\!\!\!\!\!\!\!\!
\sum_{\ \ \sigma\in S_n} \!\!\!\!\! \epsilon(\sigma\!) \!\!\!
\sum_{(\!k_i,a_i\!)}\!\!\!
\sum_{\ \ \pi\in S_r} \!\!\!\! \epsilon(\pi)(\!-\!1)^{\alpha_1}\!
f\{1^{k_0}\!,g_{\sigma(1)}^{a_\sigma(1)}\!,1^{k_1}\!,\ttdots,
g_{\sigma(n)}^{a_\sigma(n)}\!,\!1^{k_n\!}\}
\big(
x_{\pi(1+A)},\ttdots,x_{\pi(k_0+A)},
X_{\!\sigma(1)\!},
\\[-6pt]
\phantom{******************************,}
x_{\pi(1+k_0+A)},
\ttdots,
X_{\!\sigma(n)},
\\[2pt]
\phantom{**********************************,}
x_{\!\pi(1+ k_0 + \dots + k_{n\!}+A)},
\ttdots,
x_{\!\pi(r)\!}
\big)\!,\!\!\!\!\!
\end{multline*}
where $\alpha_1$ is given in Lemma \eqref{TechLemPiHat}.
Finally, when the $x$-terms are moved inside, the left hand side is equal to

\begin{multline*}
\!\!\!\!\!\!\!\!\!
\sum_{\ \ \sigma\in S_n} \!\!\! \epsilon(\sigma\!) \!
\sum_{(\!k_i,a_i\!)}\!\!\!
\sum_{\ \ \pi\in S_r} \!\!\!\! \epsilon(\pi) (\!-\!1)^{\alpha_1+\mu}
f\big\{
x_{\pi(1+A)},\ttdots,x_{\pi(k_0+A)},
g_{\sigma(1)}(X_{\sigma(1)}),
x_{\pi(1+k_0+A)},
\ttdots\\[-6pt]
\phantom{************************}
\ttdots,
g_{\sigma(n)}(X_{\sigma(n)}),
x_{\pi(1+ k_0 + \dots + k_{n}+A)},
\ttdots,
x_{\pi(r)}
\big\}\!.\!\!\!\!\!
\\[3pt]
\text{Here, }
\mu = \sum_i g_{\sigma(i)} [x_{\pi(1+A)} + \dots + x_{\pi(k_0 + \dots + k_{i-1} + A}]
+ \sum_{j<i} g_{\sigma(i)} \vert X_{\sigma(j)}\vert
\phantom{*********,}
\\
\text{and }
\epsilon(\sigma) = (-1)^{\zeta},\ \text{where }
\zeta
=
\!\!\!\!\!\!\!\!\!\!\!\!\!\!\!\!\!\!\!\!\!\!\!\!\!\!\!\!\!
\sum_{\qquad\qquad\ \ \ i<j\  \&\ \sigma(i) > \sigma(j)} 
\!\!\!\!\!\!\!\!\!\!\!\!\!\!\!\!\!\!\!\!\!\!\!\!\!\!
g_{\sigma(i)} g_{\sigma(j)}. 
\phantom{**********************}
\end{multline*}
Now that the terms on both sides are easy to compare, it is clear that
the two sides are equal if and only if 
$\tilde{\nu}+\lambda+\eta+\zeta+\alpha_1+\mu\equiv 0\ (mod\ 2)$.

After making the most obvious cancellations, we see that
$\tilde{\nu}+\lambda+\eta+\zeta+\alpha_1+\mu$
is congruent to
$$
\sum_{j<i} g_i \vert X_j\vert
+
\!\!\!\!\!\!\!\!\!\!\!\!\!\!\!\!\!\!\!\!\!\!\!\!\!\!\!\!\!
\sum_{\qquad\qquad\ \ \ i<j\  \&\ \sigma(i) > \sigma(j)} 
\!\!\!\!\!\!\!\!\!\!\!\!\!\!\!\!\!\!\!\!\!\!\!\!\!\!
\left(g_{\sigma(i)} \vert X_{\sigma(j)}\vert
+ g_{\sigma(j)} \vert X_{\sigma(i)}\vert\right)
+ \sum_{j<i} g_{\sigma(i)} \vert X_{\sigma(j)}\vert,
$$
which is congruent to zero $(mod\ 2)$ by Lemma \eqref{techlem3}.
\end{proof}

\medskip

\section{Symmetrization of the Brace Structure on $\bigoplus_{k\geq 1}Hom(V^{\otimes k},V)$}\label{secTheorem2}

In this section, we will demonstrate a nice
relationship between the the brace defined in 
Example \ref{exbrace} and the symmetric brace defined in 
Example \ref{exsymbrace}, by showing that
the symmetrization of the non symmetric brace 
structure on $Hom(V^{\otimes k},V)$ is equal to the symmetric brace of 
the anti-symmetrized maps.
Specifically, we have
\begin{theorem}\label{Theorem2}
$
\sum_{\sigma\in S_n} \epsilon(\sigma) as(f\{g_{\sigma(1)},\dots,g_{\sigma(n)}\})
=
as(f)\big<as(g_1),\dots,as(g_n)\big>.
$
\end{theorem}

\begin{proof}
First, we will manipulate the right hand side.  Using the 
symmetric brace structure defined in Example \ref{exsymbrace},
$as(f)\big<as(g_1),\dots,as(g_n)\big>(x_1,\dots,x_r)$
is equal to
$$
(-1)^{\delta}
\!\!\!\!\!\!\!\!\!\!\!
\sum_{\substack{\gamma \text{ is an }\\ 
(a_1\vert a_2\vert\dots\vert a_{n+1}) \\ 
\text{ unshuffle}}}\!\!\!\!\!\!\!\!\!\!\!\!
\chi(\gamma)
as(f)(as(g_1)\otimes\dots\otimes as(g_n)\otimes 1^{\otimes N-n})
(x_{\gamma(1)},\dots,x_{\gamma(r)}),
$$
where $\delta$ is given in Example \ref{exsymbrace}.

When we substitute the $x$ terms using the Koszul convention
and suppress the tensor notation, this
is equal to 
\begin{multline*}
\qquad\qquad (-1)^{\delta}
\!
\sum_{\gamma}
\chi(\gamma)
(-1)^{\nu}
as(f)(h_1,\ttdots, h_n,
x_{\gamma(1+\sum_{i=1}^{n} a_i))},\ttdots,x_{\gamma(r)}),
\\
\text{where }
\nu = \sum_{i=2}^n q_i [x_{\gamma(1)}+\dots + x_{\gamma(a_1 +\dots +a_{i-1})}]
\phantom{***********************}
\\[2pt]
\text{and }
h_k = as(g_k)\left(x_{\gamma\left(1+a_1 + \dots + a_{k-1}\right)},\ttdots, 
x_{\gamma\left(a_1 + \dots + a_{k}\right)}\right)
\phantom{*******************}
\\[6pt]
= 
\sum_{\ \pi_k\in S_{a_k}}\chi(\pi_k)
g_k\left(x_{\gamma\left(\pi_k(1)+a_1 + \dots + a_{k-1}\right)},\ttdots, 
x_{\gamma\left(\pi_k(a_k) + a_1 + \dots + a_{k-1}\right)}\right).\qquad\ \,
\end{multline*}

If we denote $A=\sum_{i=1}^{n} a_i$ and apply Lemma \eqref{techlem1}, 
this is equal to
$$
\!
\sum_{\gamma}
\chi(\gamma)
(-1)^{\delta + \nu\!\!}
f\circ\Phi_{\!na}\circ\Psi_{\!n}\!\!\left(
\!\!\!\!\!\!\!\!\!\!\!\!\!\!\!\!\!\!\!\!\!\!\!\!\!\!
\sum_{\qquad\qquad\ \pi_{a_n+1}\in S_{a_{n+1}}}
\!\!\!\!\!\!\!\!\!\!\!\!\!\!\!\!\!\!\!\!\!\!\!\!\!
\chi(\pi_{n+1\!}) 
(h_{1\!},\ttdots, h_{n\!},
x_{\gamma(\pi_{a_n+1}(1\!)+A\!)},\ttdots,x_{\gamma(\pi_{a_n+1}(a_n+1\!)+A\!)}
\!\!\right)\!\!.
$$
Now, we will pull all of the $x$ terms back out, in order to apply 
Lemma \eqref{techlem2}.  Note that the Koszul signs from this 
transformation merely cancel out $(-1)^{\nu}$.
We then have the following long formula, which spans two lines!
\vspace{2pt}
\begin{multline*}
(-1)^{\delta}
\sum_{\gamma}
\chi(\gamma)
\!\!\!\!\!\!\!
\sum_{\quad\ \pi_{1}\in S_{a_1}}
\!\!\!\!\!\!\!\!\!
\chi(\pi_{1}) 
\dots\!\!\!\!\!\!\!\!\!\!\!\!\!\!\!\!\!\!\!\!\!\!\!\!\!\!
\sum_{\qquad\qquad\ \pi_{(a_n+1)}\in S_{a_{n+1}}}
\!\!\!\!\!\!\!\!\!\!\!\!\!\!\!\!\!\!\!\!\!\!\!\!\!
\chi(\pi_{n+1}) 
f\circ\Phi_{\!na}\circ\Psi_n(g_1,\ttdots,g_n,1^{a_{n+1}})
\\[2pt]
\left(\!
x_{\!\gamma\left(\pi_{\!1}(1)\right)},\ttdots, 
x_{\!\gamma\left(\pi_{\!1}(a_1)\right)},
x_{\!\gamma\left(\pi_2(1)+a_1\right)},\ttdots, 
x_{\!\gamma(\pi_{a_n\!}(A\!)},
x_{\!\gamma(\pi_{a_n\!+\!1}(1)+A\!)},\ttdots,
x_{\!\gamma(\pi_{a_n\!+\!1}(a_n+1)+A\!)}
\right)\!\!.\\[-8pt]
\end{multline*}

Now, though, we can apply Lemma \eqref{techlem2},
which yields the much shorter formula,
$$
(-1)^{\delta}
\sum_{\pi\in S_{r}}
\!
\chi(\pi) 
f\circ\Phi_{\!na}\circ\Psi_{\!n}(g_1,\ttdots,g_n,1^{a_{n+1}})
\left(\!
x_{\!\pi(1)},\ttdots, 
x_{\!\pi(r)}
\right)\!.
$$

Before continuing, we need to pull all of the $x$ terms back inside.
In order to make our expressions a bit shorter,
let $X_i$ denote the input to $g_i$, and let $X_{n+1}$ denote
the free $x$ terms (letting $a_{n+1}=N\!-\!n$).  In other words, define
$$
X_i = x_{\pi\left(1+a_1 + \dots + a_{i-1}\right)},\dots, 
x_{\pi\left(a_1 + \dots + a_{i}\right)}.
$$
It will also be convenient to let $\vert X_i\vert$ denote the sum
of the degrees of the variables in $X_i$.
When we pull the $x$-terms inside and use the more concise notation 
just defined, the formula for the right hand side becomes
$$
(-1)^{\delta}\!
\sum_{\pi\in S_{r}}
\!
\chi(\pi) (-1)^{\tilde{\nu}}
f\circ\Phi_{\!n,N\!-\!n}\circ\Psi_{\!n}
\left(
g_{1}(X_1),\ttdots, g_{n}(X_n), X_{n+1}
\right)\!,
$$
where $\tilde{\nu} = \sum_{j<i} q_i \vert X_j\vert$.
After expanding $\Psi_n$,
the right hand side is equal to
$$
(-1)^{\delta}\!
\sum_{\pi\in S_{r}}
\!
\chi(\pi) (-1)^{\tilde{\nu}}
f\circ\Phi_{\!n,N\!-\!n}
\left(
\sum_{\ \sigma\in S_n}\chi(\sigma)
\left(
g_{\sigma(1)}(X_{\sigma(1)}),\dots, g_{\sigma(n)}(X_{\sigma(n)}),
X_{n+1}
\right)\!\!
\right)\!.
$$
In the above expression, $\chi(\sigma)$ is equal to $(-1)^{\lambda}$,
where 
$$
\lambda =
\!\!\!\!\!\!\!
\sum_{\substack{i<j\leq n,\\
\sigma(i) > \sigma(j)} }
\!\!\!
\left[
\Big( q_{\sigma(i)} + \vert X_{\sigma(i)}\vert \Big)
\Big( q_{\sigma(j)} + \vert X_{\sigma(j)}\vert \Big)
+1 \right].
$$
Now, if we expand $\Phi_{\!n,N\!-\!n}$,
we get
\begin{multline*}
\!\!\!\!\!\!\!\!\!\!\!\!\!\!\!\!\!\!\!\!\!\!\!\!\!\!
\sum_{\substack{\qquad \pi\in S_{r},\sigma\in S_n,
\\
\qquad\quad\ \  k_0+\dots +k_n = N\!-\!n
}}
\!\!\!\!\!\!\!\!\!\!\!\!\!\!\!\!\!\!\!\!\!\!
\chi(\pi) 
(-1)^{\delta + \tilde{\nu} + \lambda + \eta}
f\big(
x_{\pi(1+A)},\ttdots,x_{\pi(k_0+A)},
g_{\sigma(1)}(X_{\sigma(1)}),
x_{\pi(1+k_0+A)},
\ttdots\\[-10pt]
\qquad\qquad\qquad\qquad\qquad\qquad\ttdots,
g_{\sigma(n)}(X_{\sigma(n)}),
x_{\pi(1+ k_0 + \dots + k_{n}+A)},
\ttdots,
x_{\pi(r)}
\big),\ 
\\
\text{where }
\eta = \!\sum_{i=1}^n \!\left\{\!
\Big( q_{\sigma(i)} + \vert X_{\sigma(i)}\vert \Big)
\Big( 
x_{\pi(1+A)}+ \dots + x_{\pi(k_0 + \dots + k_{i-1}+A)}\!
\Big)
+ (n-i) k_i
\!\right\}.
\end{multline*}

Now, we will work with the left hand side of the equation.  Using the
brace defined in Example \ref{exbrace},\ \
$
\sum_{\sigma\in S_n} \epsilon(\sigma) as(f\{g_{\sigma(1)},\dots,g_{\sigma(n)}\})
(x_1,\dots,x_r)$  
is equal to
$$
\!\!\!\!
\sum_{\ \ \sigma\in S_n} \!\!\!\!\!\varepsilon(\sigma)  as\!\left(
\!\!\!\!\!\!\!\!\!\!\!\!\!\!\!\!\!\!\!\!\!\!\!\!\!\!
\sum_{\qquad\qquad\ k_0 + \dots k_n=N-n} 
\!\!\!\!\!\!\!\!\!\!\!\!\!\!\!\!\!\!\!\!\!\!\!\!\!
(-1)^{\beta}
f(1^{\otimes k_0}\!\otimes g_{\sigma(1)}\!\otimes\! 1^{\otimes k_1}\!\otimes\dots\otimes\!
1^{\otimes k_{n-1}}\!\otimes\! g_{\sigma(n)}\!\otimes\! 1^{\otimes k_n})\!\!\right)\!
(x_1,\!\tdots,x_r),
$$
where $\beta$ is given in Example \ref{exbrace}.  
Note also that the Koszul sign $\epsilon(\sigma)$ must be calculated using the degree of
$g_i$ as an element of the symmetric brace algebra (so
$\vert g_i\vert =q_i+a_i-1$). 
Thus $\epsilon(\sigma) = (-1)^{\zeta}$, where
$$
\zeta = 
\!\!\!\!\!\!\!\!\!\!\!\!\!\!\!\!\!\!\!\!\!\!\!\!\!\!\!\!\!\!
\sum_{\substack{\qquad\qquad\ \ \ i<j\  \&\ 
\sigma(i) > \sigma(j)} }
\!\!\!\!\!\!\!\!\!\!\!\!\!\!\!\!\!\!\!\!\!\!\!\!\!\!\!\!
(q_{\sigma(i)} + a_{\sigma(i)} - 1) 
(q_{\sigma(j)} + a_{\sigma(j)} - 1).
$$
If we now antisymmetrize by taking all signed permutations of the $x$'s, and 
suppress the tensor notation, this is equal to
$$
\!\!\!\!
\sum_{\ \ \sigma\in S_n} \!\!\!\!
\!\!\!\!\!\!\!\!\!\!\!\!\!\!\!\!\!\!\!\!\!
\sum_{\qquad\qquad\ k_0 + \dots k_n=N-n} 
\!\!\!\!\!\!\!\!\!\!\!\!\!\!\!\!\!\!\!\!\!\!\!\!\!
(-1)^{\beta + \zeta}
f\left(1^{k_0}, g_{\sigma(1)}, 1^{k_1},\ttdots,
1^{k_{n-1}},g_{\sigma(n)}, 1^{k_n}\right)\!
\left(\!\!\!
\sum_{\ \ \pi\in S_r} \!\!\!\!
\chi(\pi)
\left(x_{\pi(1)},\ttdots,x_{\pi(r)}\right)\!
\right)\!.
$$

After applying Lemma \eqref{TechLemPiHat}, 
the left hand side is equal to
\begin{align*}
\!\!\!\!\!\!\!\!\!\!\!\!\!\!\!\!\!\!\!\!\!\!\!\!\!\!\!\!
\sum_{\substack{\qquad\quad \pi\in S_{r},\sigma\in S_n,
\\
\qquad\qquad\ k_0+\dots +k_n = N\!-\!n
}}
\!\!\!\!\!\!\!\!\!\!\!\!\!\!\!\!\!\!\!\!\!\!\!\!\!
(-1)^{\beta + \zeta+\alpha_2}\chi(\pi)
f\!\left(1^{k_0}, g_{\sigma(1)}, 1^{k_1},\ttdots,
g_{\sigma(n)}, 1^{k_n}\!\right)\!
\big(
&x_{\pi(1+A)},\ttdots,x_{\pi(k_0+A)},
X_{\sigma(1)},\\[-14pt]
&x_{\pi(1+k_0+A)},
\ttdots,
X_{\sigma(n)},\\
&x_{\pi(1+ k_0 + \dots + k_{n}+A)},
\ttdots,
x_{\pi(r)}
\big)\!,
\end{align*}
where $\alpha_2$ is given in Lemma \eqref{TechLemPiHat}.

Finally, when the variables are moved inside, the left hand side is 
equal to 
\begin{multline*}
\!\!\!\!\!\!\!\!\!\!\!\!\!\!\!\!\!\!\!\!\!\!\!\!\!\!\!\!
\sum_{\substack{\qquad\quad \pi\in S_{r},\sigma\in S_n,
\\
\qquad\qquad\ k_0+\dots +k_n = N\!-\!n
}}
\!\!\!\!\!\!\!\!\!\!\!\!\!\!\!\!\!\!\!\!\!\!\!\!\!
(-1)^{\beta + \zeta+\alpha+\mu}\chi(\pi)
f\big(
x_{\pi(1+A)},\ttdots,x_{\pi(k_0+A)},
g_{\sigma(1)}(X_{\sigma(1)}),
x_{\pi(1+k_0+A)},
\ttdots\\[-10pt]
\qquad\qquad\qquad\qquad\qquad\quad\ttdots,
g_{\sigma(n)}(X_{\sigma(n)}),
x_{\pi(1+ k_0 + \dots + k_{n}+A)},
\ttdots,
x_{\pi(r)}
\big),\ 
\\[8pt]
\text{where }
\mu = \sum_i q_{\sigma(i)} [x_{\pi(1+A)} + \dots + x_{\pi(k_0 + \dots + k_{i-1} + A}]
+ \sum_{j<i} q_{\sigma(i)} \vert X_{\sigma(j)}\vert.
\qquad\qquad
\end{multline*}

Since the right hand side is equal to 
\begin{multline*}
\!\!\!\!\!\!\!\!\!\!\!\!\!\!\!\!\!\!\!\!\!\!\!\!\!\!
\sum_{\substack{\qquad \pi\in S_{r},\sigma\in S_n,
\\
\qquad\quad\ \  k_0+\dots +k_n = N\!-\!n
}}
\!\!\!\!\!\!\!\!\!\!\!\!\!\!\!\!\!\!\!\!\!\!
\chi(\pi) 
(-1)^{\delta + \tilde{\nu} + \lambda + \eta}
f\big(
x_{\pi(1+A)},\ttdots,x_{\pi(k_0+A)},
g_{\sigma(1)}(X_{\sigma(1)}),
x_{\pi(1+k_0+A)},
\ttdots\\[-10pt]
\ttdots,
g_{\sigma(n)}(X_{\sigma(n)}),
x_{\pi(1+ k_0 + \dots + k_{n}+A)},
\ttdots,
x_{\pi(r)}
\big),\ 
\end{multline*}
we see that the two sides are equal if and  only if 
$$\beta + \zeta+\alpha_2+\mu + \delta + \tilde{\nu} + \lambda + \eta \equiv 0\ (mod\ 2).$$

After cancelling the most obvious terms,
$\beta\!+\!\zeta\!+\!\alpha_2\!+\!\mu\!+\!\delta\!+\!\tilde{\nu}\!+\!\lambda\!+\!\eta$
is congruent to

\begin{multline*}
\sum_i (n-i) a_{\sigma(i)}
+ \sum_i(N-i)q_{\sigma(i)} +\sum_{j < i} q_{\sigma(i)} a_{\sigma(j)}
\\+
\!\!\!\!\!\!\!\!\!\!\!\!\!\!\!\!\!\!\!\!\!\!\!\!\!\!\!\!\!\!
\sum_{\substack{\qquad\qquad\ \ \ i<j\  \&\ 
\sigma(i) > \sigma(j)} }
\!\!\!\!\!\!\!\!\!\!\!\!\!\!\!\!\!\!\!\!\!\!\!\!\!\!\!\!
[
q_{\sigma(i)} a_{\sigma(j)} + q_{\sigma(i)}
+ a_{\sigma(i)}q_{\sigma(j)} + a_{\sigma(i)}
+ q_{\sigma(j)} + a_{\sigma(j)}
]
+
\sum_{j<i} q_{\sigma(i)} \vert X_{\sigma(j)}\vert
\\+
\sum_i (N-i) q_i 
+
\sum_{j<i}q_i a_j + \sum_i (n-i) a_i
+
\sum_{j<i} q_i \vert X_j\vert
+
\!\!\!\!\!\!\!\!\!\!\!\!\!\!\!\!\!\!\!\!\!\!\!\!\!\!\!\!\!\!
\sum_{\substack{\qquad\qquad\ \ \ i<j\  \&\ 
\sigma(i) > \sigma(j)} }
\!\!\!\!\!\!\!\!\!\!\!\!\!\!\!\!\!\!\!\!\!\!\!\!\!\!\!\!
\left[q_{\sigma(i)}\vert X_{\sigma(j)}\vert
+ \vert X_{\sigma(i)}\vert q_{\sigma(j)}
\right].
\end{multline*}

After applying Lemma \eqref{techlem3}, this is congruent to

$$
\sum_i 
\left\{(n\!-\!i) a_{\sigma\!(\!i\!)} + (N\!-\!i)q_{\sigma\!(\!i\!)} 
+
(i\!-\!1)\left[ a_i + a_{\sigma\!(\!i\!)}+q_i + q_{\sigma\!(\!i\!)} \right]
+ (N\!-\!i) q_i + (n\!-\!i) a_i
\right\},
$$

which is equal to
$\sum_i \left\{
(n-1)\left[ a_{\sigma(i)} + a_i \right]
+ (N-1)\left[ q_{\sigma(i)} + q_i \right]
\right\}\equiv 0\ (mod\ 2)$.
\end{proof}

As a corollary, we obtain Theorem 3.1 of \cite{LM95}:
\begin{corollary}
The anti-symmetrization $l:=as(\mu)$ of an $A_{\infty}$- algebra structure $\mu$ yields an $L_{\infty}$-algebra structure.
\end{corollary}
\begin{proof} Given $\mu\{\mu\}=0$ (recall Remarks \eqref{Ainf} and\eqref{Linf}) we have
$$0=as(\mu\{\mu\})=as(\mu)\langle as(\mu)\rangle=l\langle l\rangle.$$
\end{proof}

\section{Acknowledgements}

We would like to thank  Martin Markl for  providing many helpful suggestions. The first author would also like
to thank the Mathematics Institute of the Czech Academy for its hospitality during her visit to
Prague in September, 2003.

\end{document}